\documentclass{amsart}

\usepackage{amsmath}
\usepackage{amssymb}
\usepackage{a4wide}
\usepackage{graphicx}

\renewcommand{\rho}{\varrho}

\renewcommand{\phi}{\varphi}

\begin{document}

\title{Ten colours in quasiperiodic and regular hyperbolic tilings}

\author{R. L\"uck}
\address{Stuttgart, Germany}
\email{r.v.lueck@web.de}

\author{D. Frettl\"oh}
\address{Fakult\"at f\"ur Mathematik, Universit\"at Bielefeld,
Postfach 100131, 33501 Bielefeld,  Germany}
\email{dirk.frettloeh@math.uni-bielefeld.de}
\urladdr{http://www.math.uni-bielefeld.de/baake/frettloe}

\begin{abstract}
Colour symmetries with ten colours are presented for different tilings. 
In many cases, the existence of these colourings were predicted by
group theoretical methods. Only in a few cases explicit constructions
were known, sometimes using combination of two-colour and
five-colour symmetries. Here we present explicit constructions of 
several of the predicted colourings for the first time, and discuss
them in contrast to already known colourings with ten colours.
\end{abstract}


\maketitle


\section{Introduction} 

The discovery of quasicrystals generated interest in colour symmetries
of quasicrystallographic structures, see \cite{baa}, \cite{luc} and
references therein. 
In \cite{baa} and \cite{bg}, the existence and the number of
colour symmetries using $k$ colours is determined for a large class of
quasicrystallographic structures. Despite of the richness of regular
structures in hyperbolic plane, colour symmetries of  
regular tilings in the hyperbolic plane are not widely studied. For
an exception see \cite{mlp} and references therein. Recently a method
was presented to determine the possible number of colours for which a
colour symmetry of any regular tiling exists, either in the Euclidean,
spherical or hyperbolic plane \cite{fl}. This purely
algebraic method does not yield the colouring itself. 
Here we focus on colour symmetries with ten colours. We present
explicit colourings with ten colours which are predicted by
\cite{fl}. 

Basic definitions about regular tilings are quite standard and can be 
found in \cite{gs}, for instance, from where we also adopt the notation. 
A tiling of the plane is said to be {\em regular}, if all tiles
are regular polygons, and all vertex figures are regular polygons.
Any ordered pair of integers $(p,q)$, where $p,q \ge 3$, defines a
regular tiling by regular $p$-gons, where $q$ tiles meet at 
each vertex. Following \cite{gs}, such a tiling is denoted by
$(p^q)$. The value $d:= 2(\frac{1}{p} + \frac{1}{q})$ determines in
which plane $(p^q)$ lives: If $d>1$ then $(p^q)$ is a regular
spherical tiling (which can be regarded as a Platonic solid). If
$d=1$ then $(p^q)$ is a regular Euclidean tiling, and if 
$d<1$ it is a regular tiling of the hyperbolic plane. 
Figures in this paper represent regular tilings of the hyperbolic 
plane, using the Poincar\'e disc model.

\section{Colour symmetries of regular tilings with ten colours}

A symmetry of a tiling is any isometry which maps the tiling to
itself. In the sequel we distinguish two cases: Either we
consider all isometries (reflections, rotations, and their products),
or direct isometries only (no reflections, only rotations and their
products). By colour symmetry we denote here any colouring of
a tiling such that any symmetry of the uncoloured tiling acts as a
global permutation of colours. Such a colour symmetry is called
'perfect colouring' in \cite{gs}. In Table \ref{tabelle} we list the
regular tilings which are known to possess colour symmetries with ten
colours, first with respect to all isometries, then with respect to
direct isometries only.

\begin{table}[h]
{\small 
\begin{tabular}{l|l|l|lllllllllllll}
          & spher & eucl & hyp \\
\hline
all isometries & $(3^5)$ &  & $(3^8)_2$  & & $(3^{10})_3$  & $(4^5)_4$
  & $(4^6)$ & $(6^4)_2$ & $(6^5)_6$ & & & 
\\ 
direct isom.  &  $(3^5)$ & $(4^4)$ & $(3^8)_3$  & $(3^9)_3$  &
$(3^{10})_6$  & $(4^5)_6$  & $(4^6)_3$ & $(6^4)_6$ & $(6^5)_{15}$ &
$(4^9)_8$ & $(8^3)$ & $(8^4)_7$ 
\\ 
\hline 
\end{tabular} }
\vspace{1mm}
\caption{Regular tilings which have a 10-colour symmetry. A lower
  index denotes the number of different colour symmetries for that
  tiling, where no lower index means 'one'. \label{tabelle}}
\end{table}

\begin{figure}
\begin{center}
\includegraphics[width=100mm]{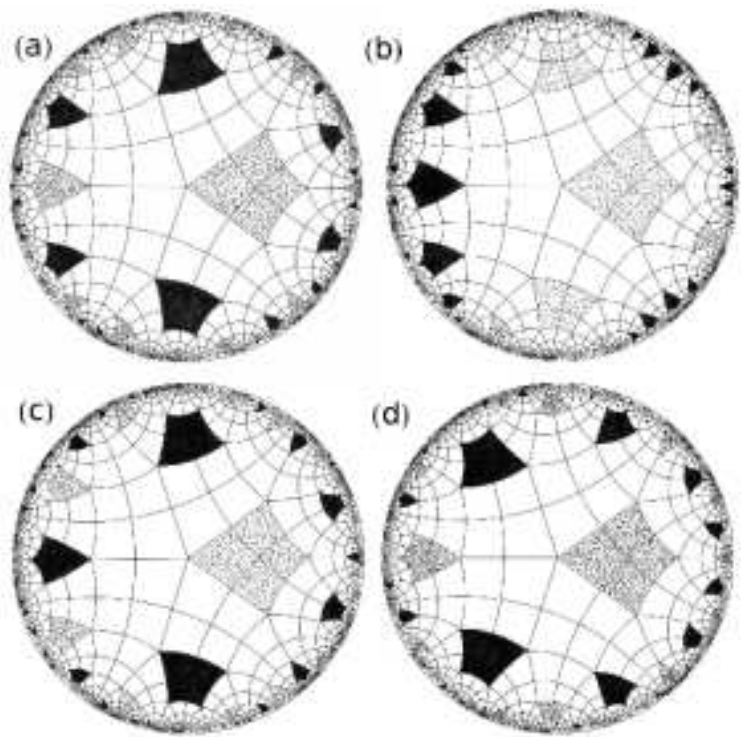}
\end{center}
\caption{Two out of ten colours of a perfect colour symmetry in the
hyperbolic tiling $(4^5)$, indicated by black and gray. Identifying 
black and gray yields a five-colour symmetry. 
(a) - (c) are based on the same five-colour symmetry, 
(d) is based on a second five-colour symmetry. See Fig.\ \ref{a1} for
a colour version.
\label{eins}} 
\end{figure}

\begin{figure}
\begin{center}
\includegraphics[width=55mm]{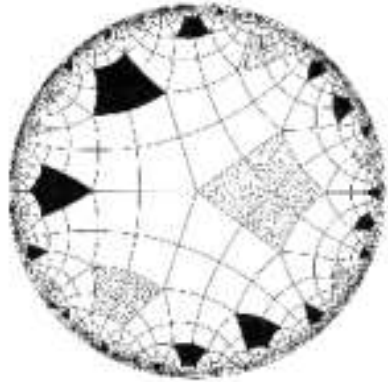}
\end{center}
\caption{Example of a 10-colour symmetry in the hyperbolic tiling 
$(4^5)$ which is only perfect with respect to rotational symmetry. 
Only two out of ten colours are marked. The combination of black and
gray results in the five-colour symmetry of the second type.
See Fig.\ \ref{a2} for a colour version.
\label{zwei}} 
\end{figure}

\begin{figure}
\begin{center}
\includegraphics[width=140mm]{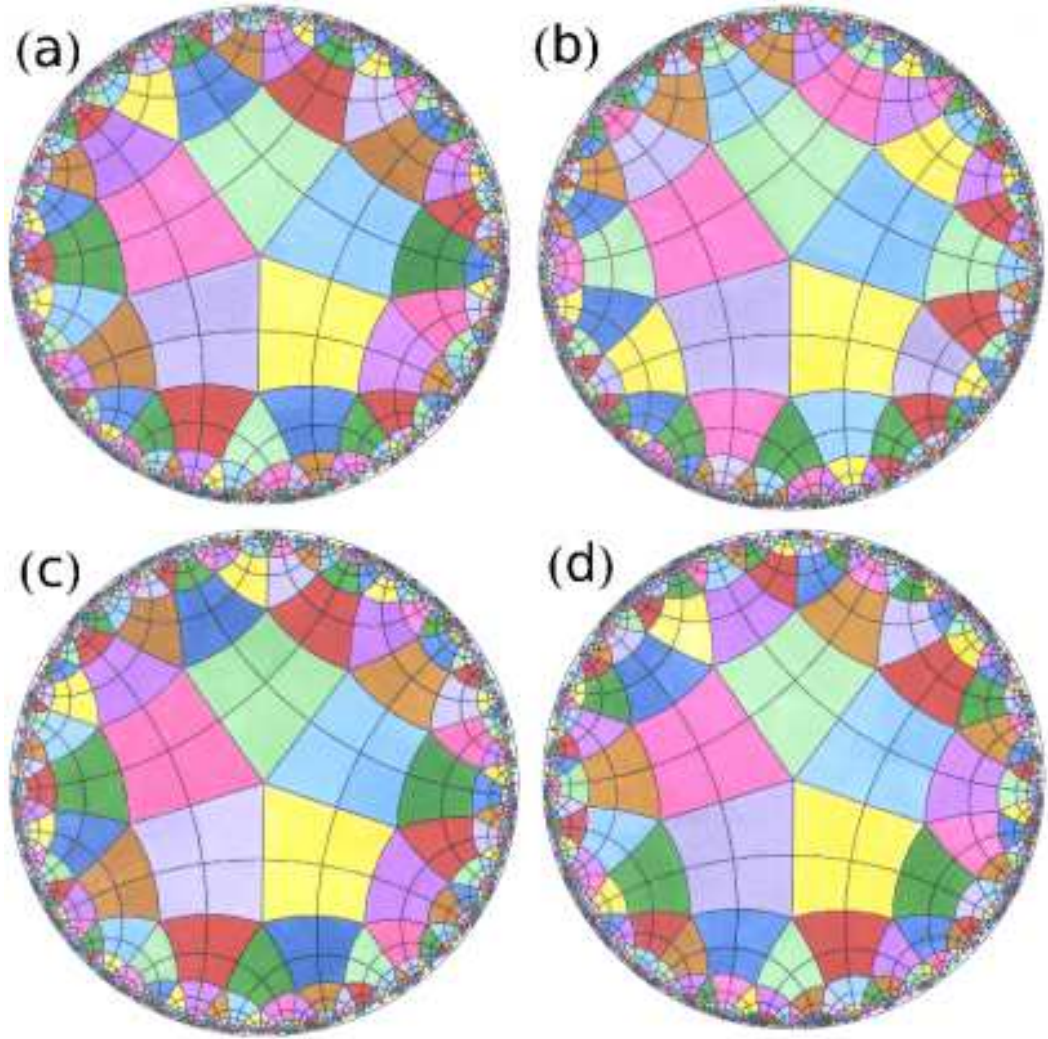}
\end{center}
\caption{Perfect colour symmetry in the hyperbolic tiling $(4^5)$
  designated by 10 colours. The combinations of dark and light colours
  represent a five-colour symmetry. (a) - (c) are based on a first
  type of a five-colour decoration, (d) is based on a second type of a
  five-colour decoration (compare Fig.\ \ref{eins}).    
\label{a1}} 
\end{figure}

\begin{figure}
\begin{center}
\includegraphics[width=140mm]{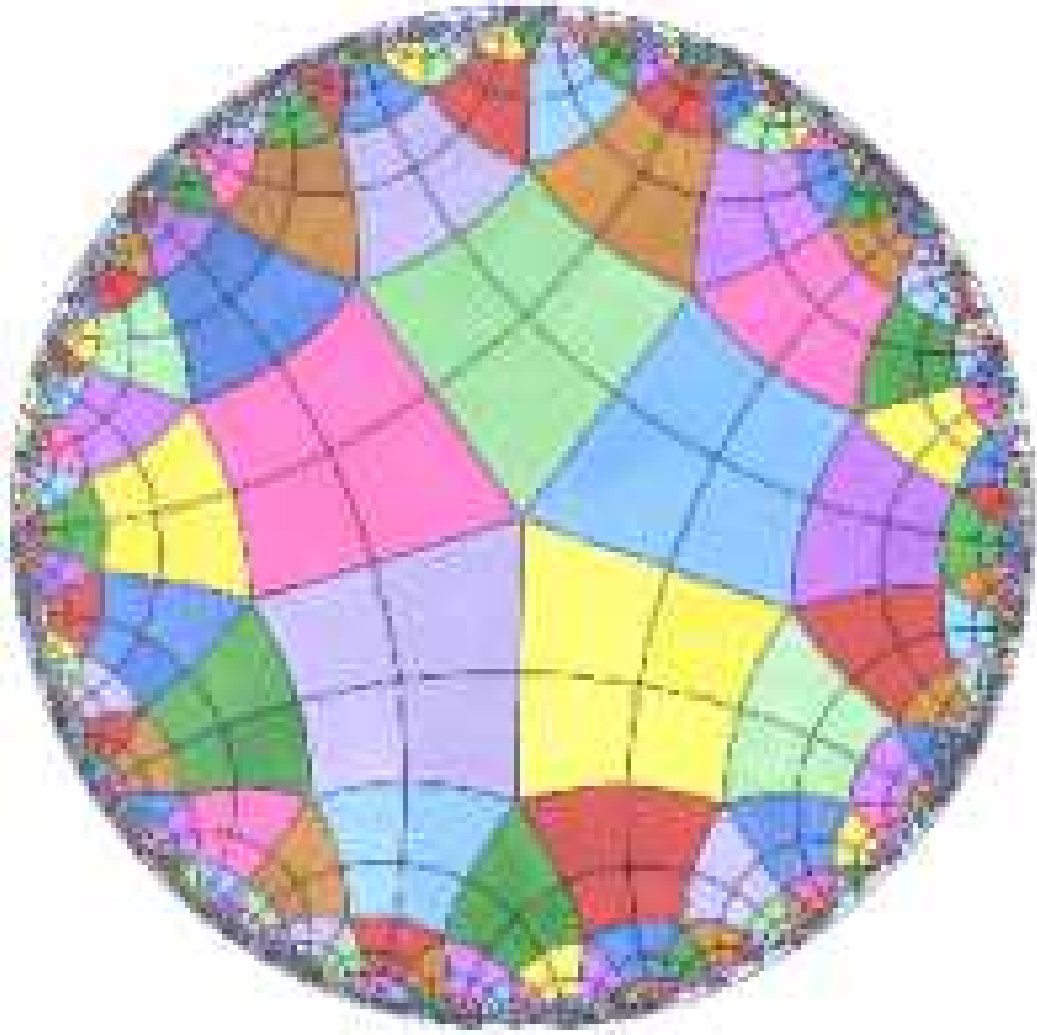}
\end{center}
\caption{Example of a 10-colour symmetry in the
  hyperbolic tiling $(4^5)$ which is only perfect with respect to
  rotational symmetry. The combination of dark and light colours
  results in the five-colour symmetry of the second type (compare
  Fig.\ \ref{zwei}). 
\label{a2}}
\end{figure}

There is only one 10-colour symmetry for regular spherical tilings:
for $(3^5)$, which can be regarded as a coloured
icosahedron. Moreover, there is only one 10-colour symmetry for 
regular euclidean tilings: 
for the square tiling $(4^4)$. 
In contrast, there are many 10-colour symmetries with ten colours for
regular hyperbolic tilings. For instance, there are four of those for
the tiling $(4^5)$ with respect to all isometries. These are shown in
Figure \ref{eins}. (Note, that all figures show the tiling $(p^q)$
superimposed with its dual tiling $(q^p)$).   
Two of ten colours are marked by black and gray, the other colours 
are omitted for the sake of clarity and are resulting from fivefold 
rotation around the centre of the figure. Combining black 
and gray yields a five-colour symmetry. Figures (a) - (c) are based 
on the same five-colour symmetry. Figure (d) is based on a
second five-colour symmetry. A full-colour version of Fig. \ref{eins}  
is given in Fig.\ \ref{a1} of the on-line version of the present paper.
For the same tiling, there are six colour symmetries with respect to
direct isometries, including the former four. An additional one is
shown in Figure \ref{zwei}. A full colour version is represented in
Fig.\ \ref{a2} in the online version. It is based on the second
five-colour symmetry mentioned above. The figure represents one
partner of an enantiomorphic pair. Note, that a reflection in the
horizontal axis does not permute entire colour classes, which shows
that this one is not a colour symmetry with respect to all isometries. 



\begin{figure}
\begin{center}
\includegraphics[width=100mm]{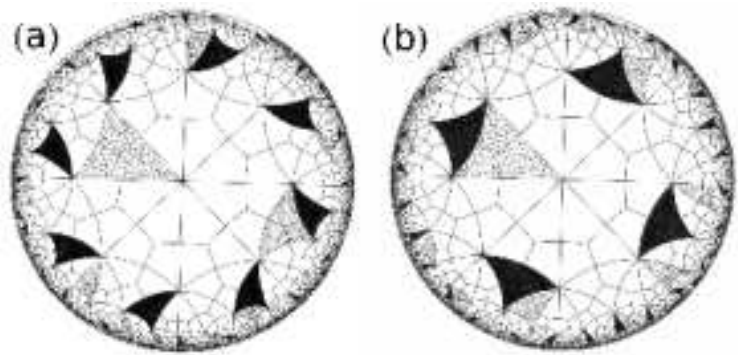}
\end{center}
\caption{Perfect colour symmetry with ten colours in the hyperbolic
tiling $(3^8)$. The gray colour symbolyses the symmetry of eight
colours in the orbit around the eightfold centre in both figures (a)
and (b). The black colour in Fig.\ (a) stands for a single colour in
two orbits around the eightfold centre, in Fig.\ (b) it stands for two
colours in a single orbit. 
\label{drei}} 
\end{figure}

According to Table \ref{tabelle}, there are two different perfect
colourings with  
ten colours for the regular tiling $(3^8)$. Figure \ref{drei} represents 
both examples for an eightfold centre. The gray colour represents 
one of eight colours in the orbit around the eightfold centre in both 
\ref{drei}(a) and \ref{drei}(b). The black colour in \ref{drei}(a)
stands for a single colour in two orbits around the eightfold centre,
in \ref{drei}(b) it stands for two colours in a single orbit.

\begin{figure}
\begin{center}
\includegraphics[width=55mm]{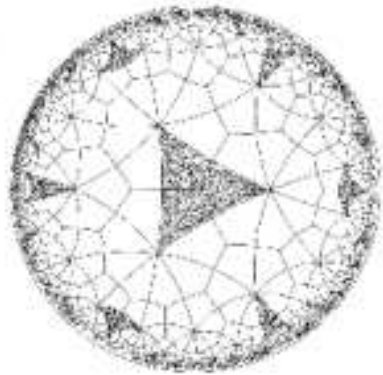}
\end{center}
\caption{A single colour of Fig.\ \ref{drei} a represented for a 
threefold centre. This figure results also as a coincidence site 
lattice of threefold vertices for $60^{\circ}$ or $180^{\circ}$
rotations around the figure centre.    
\label{vier}} 
\end{figure}

Figure \ref{drei}(a) is based on a coincidence site lattice (CSL, see
\cite{bghz}). This can be seen if the threefold centre is
transformed to the centre of the figure (see Figure \ref{vier}) and if
the pattern is rotated by $60^{\circ}$ or $180^{\circ}$: points which
are vertices in both the rotated and the unrotated image are elements
of the CSL. In Figure \ref{vier} the triangles surrounding a threefold
coinciding vertex are coloured.

\section{Colour symmetries of quasiperiodic tilings with ten colours}

Ten colours in $n$-fold symmetric Euclidean tilings are predicted only
for the periodic square lattice \cite{baa}, \cite{bg}. In several
10-fold tilings, including Penrose tilings, synthetic combinations of
five-colour symmetry with two-colour symmetry is possible. The colour
symmetry with five colours has been discussed several times
(\cite{luc} and references therein), a two-colour symmetry for
decagonal tilings was presented by Li et al.\ \cite{ldk}. The
combination of both colour symmetries in decagonal tilings is not
straightforward. Both colour symmetries are restricted to special
cases \cite{luc} and often these cases exclude each other. A possible
combination is represented by the 'orientated points' in the Penrose
pentagon pattern as published by Gummelt \cite{gum}. 

\section{Concluding remarks}

Colour symmetries of regular Euclidean patterns are composed of
operations  based on prime powers \cite{baa} and their sequence can
usually be permuted. (Usually means, that the prime does not ramify
over the underlying integer ring of $n$-fold cyclotomic integers.) In
particular, if there are $k$ colourings with $p$ colours and $m$
colourings with $q$ colours, where $p$ and $q$ are distinct primes,
then there are $km$ colour symmetries with $pq$ colours.
The composition of colour symmetries in hyperbolic space can be quite
different: a genuine colour symmetry is not necessarily based on a
prime power. Composition of several numbers of colours is often found.
In contrast to the Euclidean case, the operations can not be permuted
in general. However, in some cases they can. For instance, 
the decoration of $(4^6)$ with all isometries according to the
table can be composed from a two-colour and a five-colour symmetry.

The role of coincidence site lattices is also
different. In planar Euclidean space the CSLs exist as enantiomorphic
pairs  and force a colouring which is only perfect with respect to
rotations.  For the hyperbolic space, we found enantiomorphic pairs
which are not based on CSLs, and CSLs which are not related to a
colour symmetry.

\end{document}